\documentclass[12pt]{article}
\input epsf
\usepackage{epsfig}
\def\R{\mathbf{R}}
\def\C{\mathbf{C}}
\title{Simultaneous stabilization, avoidance and Goldberg's constants}
\author{Alex Eremenko\thanks{Supported by NSF grant DMS1067886.}}
\begin{document}
\maketitle
\begin{center}
Dedicated to the memory of A. Goldberg and V. Logvinenko
\vspace{12pt}
\end{center}
\begin{abstract}This is an exposition for mathematicians
of some unsolved problems arising in control theory of
linear time-independent systems.
\end{abstract}

The earliest automatic control devices that I know are
described in the book of Hero of Alexandria ``Pneumatica'', see
Fig. 1. In the modern times these devices are omnipresent
(almost every home appliance contains at least one, a car has several,
an airplane or a guided missile has many;
an ingenious mechanical steering device of a sailboat permits you
to sleep and to dine during your voyage, while it keeps prescribed
direction with respect to the wind; one can add many other examples).

The mathematical theory of these devices begins, as far as I know,
with George Biddell Airy (of the Airy function), Astronomer Royal, 
who investigated mathematically
stabilization of the clockwork mechanism directing
his equatorial.\footnote{Equatorial is a device that
continuously adjusts a telescope direction to compensate for the
diurnal rotation of the Earth.
One of the most complicated modern
control systems directs the Hubble telescope. It has to keep
the direction of the telescope with high accuracy and
for long time.}
The stability condition that ``all poles must be in the left
halfplane'' was explicitly stated for the first time by
J.~C.~Maxwell \cite{Max}.
Parallel research was done in Eastern Europe by Aurel Stodola (1894)
and 
Ivan Vyshnegradsky (1877) who pioneered the use of complex function theory,
anticipating the work of Nyquist (1932), see for example,
\cite{LS}. 

Most of the XIX century research in the area was related to
stabilizing the system which consists of a steam engine controlled
by the governor. 

These investigations led to the famous
criteria in terms of coefficients of a polynomial for all its roots
to belong to the left half-plane, (E. Routh, 1877, A. Hurwitz, 1895),
see \cite{G}.
\vspace{.2in}
\begin{center}
\epsfxsize=4.2in
\centerline{\epsffile{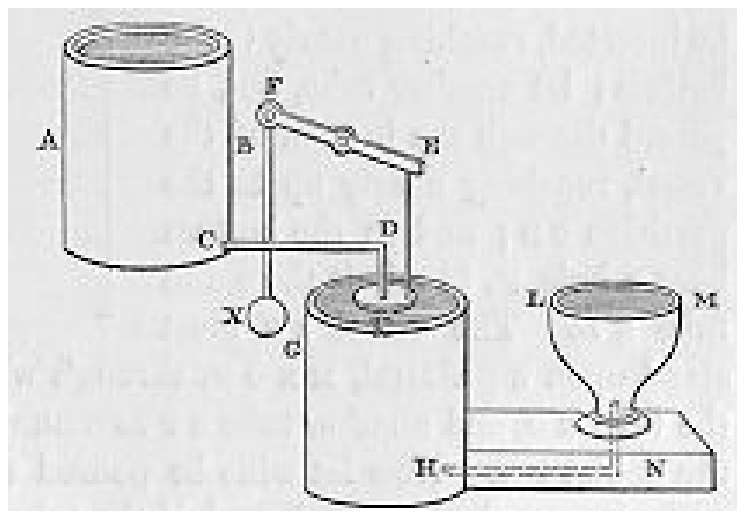}}
\vspace{.2in}\nopagebreak
Fig. 1  A XIX century illustration made according to the description
in the book of Hero Pneumatica.
\end{center}
\vspace{.2in}

A {\em linear system} of the simplest kind
is described by $3$ real matrices: $(A,B,C)$
of sizes $n\times n$, $n\times m$ and $p\times n$ respectively.
We have vectors depending on time:
the {\em inner state} $x(t)$ with values in $R^n$, the {\em input} $u(t)$
with values in $R^m$ and the {\em output} $y(t)$ with values in $R^p$.
These are related in the following way:
\begin{eqnarray*}x'&=&Ax+Bu,\\ y&=&Cx.
\end{eqnarray*}
We will only consider the case $m=p=1$
(so called single input -- single output systems). 

Taking Laplace transforms, and assuming that $x(0)=u(0)=y(0)=0$, we obtain
$zX(z)=AX(z)+BU(z),\; Y(z)=CX(z),$ so
\begin{equation}\label{1}
Y(z)=C(zI-A)^{-1}BU(z)=p(z)U(z).
\end{equation}
The rational function $p(z)=C(zI-A)^{-1}B$
is called the {\em transfer function}.
It is real and $p(\infty)=0$. Rational functions satisfying
$p(\infty)=0$ are called {\em proper}.
For every proper rational $p$ function there exists
a triple $(A,B,C)$ so that 
$p(z)=C(zI-A)^{-1}B$.

The correspondence between triples of matrices and rational functions
is not trivial, not bijective, and there is a large literature on
recovery of $A,B,C$ from the transfer function (realization theory).
But all essential properties of the system are encoded in the transfer
function and here we identify a linear system with its transfer function.

Improper transfer functions are equally important, they arise from more general
systems of differential equations with constant
coefficients; I don't go into detail,
but the primary object in this paper will
be an {\em arbitrary} real rational function; 
we call it a transfer function.
It completely describes a linear system. 

A system is called {\em stable} if the transfer function
has no poles 
in the open right half-plane $H$. The poles of the transfer function are
nothing but the eigenvalues of the matrix $A$ of the system.

For a given (maybe unstable) system, one may wish to {\em stabilize} it
by attaching a feedback {\em controller}. A controller is a linear
system of the same kind; it is described
by another real rational transfer function $c(z)$.
Attaching a controller as in the third diagram in Fig. 2 means that
we take the output of our original system, transform it by the controller,
and then add to the input:
$$Y=p(U+cY)=pU+pcY.$$
We obtain a new system, which is called the {\em closed loop system}.
By solving with respect to $Y$ we get the {\em closed loop transfer function}:
\begin{equation}\label{2}
\frac{p}{1-cp}.
\end{equation}
Cancellation between poles and zeros of $c$ and $p$ is possible
here, but engineers naturally do not want to rely on such cancellation.
So they give the following definition:

A controller $c$ {\em internally stabilizes} $p$ if $1-cp$ has no
zeros in the right half-plane $H$, the poles of $c$ are disjoint from the
zeros of $p$ in $H$, and the zeros of $c$ are disjoint from the poles of $p$
in $H$.

From now on by ``stabilization'' we mean ``internal stabilization''.
One can easily show that internal stabilization is equivalent to
the condition that all four transfer functions
$$pc/(1-pc),\quad c/(1-pc),\quad p/(1-pc),\quad 1/(1-pc)$$
are without poles in $H$.

All these four transfer functions can be realized  by attaching
the feedback in various ways, as shown in Fig. 2 below.
\vspace{.2in}

\begin{center}
\epsfxsize=4.2in
\centerline{\epsffile{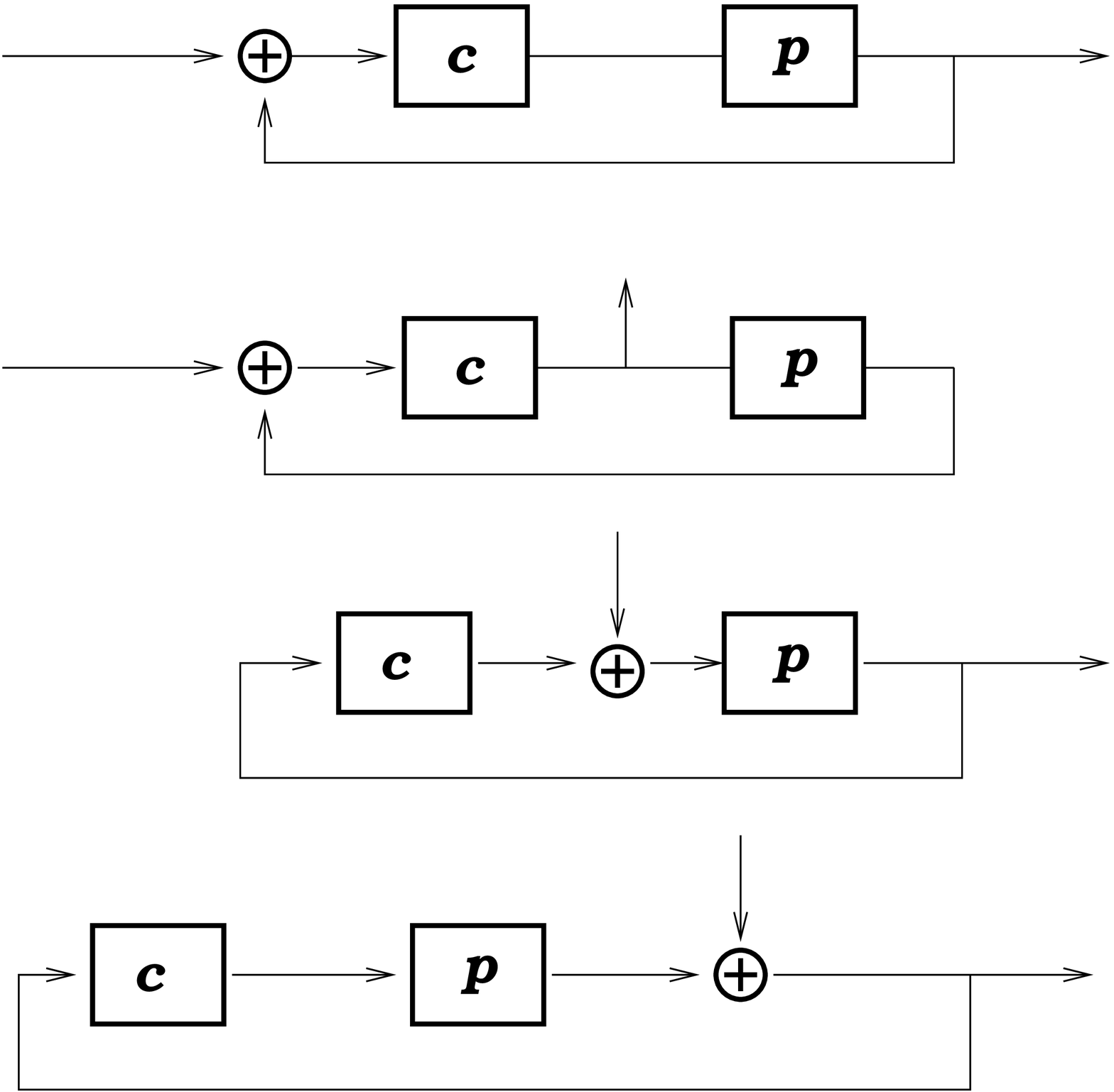}}
\vspace{.2in}

Fig. 2 Configurations corresponding to
$pc/(1-pc),\; c/(1-pc),\; p/(1-pc),\; 1/(1-pc)$.
\end{center}
\vspace{.2in}

Here is another elegant way to rewrite the internal stabilization condition:

$c$ internally stabilizes $p$ iff $c$ {\em avoids} $1/p$ in the sense
that
\begin{equation}\label{3}
c(z)\neq 1/p(z),\quad z\in H.
\end{equation}

Thus the {\em stabilization problem} is: {\em
for a given rational function $p$ find a rational
function $c$ so that (\ref{3}) holds.}

We obtain an equivalent problem when the right half-plane $H$
is replaced by the unit disc $D$.
From the point of view of system theory,
the unit disc setting corresponds
to discrete-time systems.
Instead of a differential equation, we have a
recurrence relation,
$$x(n+1)=Ax(n)+Bu(n),\quad y(n)=Cx(n),$$
and in place of the Laplace transform
we use the generating function
$X(z)=\sum_{-\infty}^\infty x(n)z^{n}$.
Then the transfer function $p(z)$ becomes $C(z^{-1}I-A)^{-1}B$,
exactly as in the case
of continuous time,
and the system is stable if $p$ has no poles in the unit disc,
which now is equivalent to saying that $A$ has no
eigenvalues whose absolute value
is greater than $1$.

Stabilization of one system is always possible if
one 
does not restrict the degree of $c$.

Now we consider simultaneous stabilization of several systems
by one controller. The problem has evident practical meaning:
the system that we want to stabilize may work in several different regimes
(think of the cooling/heating system in your home, which is usually
controlled by a single devise, or an airplane during take of/landing/horizontal
flight), mathematically this means that we want a single controller
to stabilize several systems.

Consider the problem of stabilizing two systems $p_1$ and $p_2$.
This is equivalent to
stabilization of one system by a {\em stable
controller} \cite{Blondel}.

Stability of the controller is a desirable property by itself:
if the system $p$
suddenly stops working, we don't want the controller to destroy itself.

A rational function avoiding two
different rational functions $p_i$ in $H$ 
always exists, but it may have complex coefficients,
even if $p_i$
are real. One usually needs a controller with real
coefficients. There is an obvious topological
obstruction to the existence
of real rational function without poles on $\R_{>0}$ which avoids a given
real rational function. In fact, this is the only obstruction:
\vspace{.1in}

{\em For a given real rational function $p$, there exists a real rational
function $c$ satisfying (\ref{3}) and without poles in $H$
if and only if $p$ has even number
of poles between every two adjacent zeros on $\R_{>0}$. }
\vspace{.1in}

This neat statement is due to Youla, Bongiorno and Lu
\cite{YBL}, and control theorists are very fond of it \cite{Blondel}.

Now we consider simultaneous stabilization of three systems.

Finding a function which avoids three given functions is an interesting
problem which attracted attention of pure mathematicians who were unaware of
its application to control theory.
It seems that the problem was stated for the first time
in \cite{RY}, and a connection with an ``interpolation problem''
of the kind stated below in Theorem 1 was established.

In \cite{Logv} this avoidance problem
is considered for meromorphic functions which avoid given rational
functions in an arbitrary given region. The author credits Volberg
and Eremenko who stated the problem and obtained some partial results.
Apparently they were motivated by the ``Lambda-lemma'' and
holomorphic motions which were discovered about that time \cite{MSS,L},
\cite{Slod}, \cite{Chirk}. The lambda-lemma says that if finitely
many meromorphic
functions avoid each other in a disc, that is, if their graphs
are disjoint, then one can always find an additional function which avoids all
of them.

The main conclusion in \cite{Logv}
is that in any region one can always avoid two
functions, but in general one cannot avoid three.
This is easy to explain for the case of avoidance of three
rational functions in $\C$. Take the avoided functions to be 
$0$, $\infty$ and $z$. If a meromorphic function $f$ avoids them,
then it must be rational, by Picard's Theorem, but a rational
function $f$ that avoids $0$ and $\infty$ in $\C$ must be
constant, so it cannot avoid $z$. 

Similar results were obtained in control theory for
the case of the unit disc or a halfplane.

We will work in the unit disc from now on.
First we give a general reformulation of stabilization of three
systems in terms of some unusual interpolation problem.
Various special cases of this result are mentioned in the 
control literature, but I could not find a general statement.
\vspace{.1in}

\noindent
{\bf Theorem 1.} {\em Let $\phi_1,\phi_2$ and $\phi_3$
be three rational functions without common poles,
and suppose that the set
\begin{equation}\label{three}
E=\{z\in D:\phi(z)=\phi_2(z)=\phi_3(z)\}
\end{equation}
is empty.

There exists a rational function $f$ which avoids $\phi_i$ in
$U$ if and only if there exists a rational function $g$
with the properties:

\noindent
(i)
divisor of zeros of $g$ coincides with the divisor of zeros of $\phi_3-\phi_2$;

\noindent
(ii) divisor of poles of $g$ coincides with the divisor of zeros of $\phi_3-\phi_1$,
and

\noindent
(iii)
divisor of ones of $g$ coincides with the divisor of zeros of $\phi_1-\phi_2$.
}
\vspace{.1in}

Condition that the $\phi_i$ have no common poles is added only
for simplicity of formulation: the whole situation is invariant
with respect to composition with fractional-linear transformations.
Condition that $E=\emptyset$ holds for generic $\phi_i$.

The correspondence between $f$ and $g$ is given by the
cross-ratio
$$g=\frac{(f-\phi_1)(\phi_3-\phi_2)}{(f-\phi_2)(\phi_3-\phi_1)}.$$
In the case that there are triple intersections, that is $E\neq\emptyset$,
one has to add
the condition for each point $a\in E$:
\begin{equation}\label{extra}
g(z)(\phi_3(z)-\phi_1(z))/(\phi_3(z)-\phi_2(z))=1+O(z^k),\quad z\to a,
\end{equation}
where $k$ is the order of the zero of $\phi_1-\phi_2$ at $a$.

Thus, simultaneous stabilization of three functions
(and the problem of avoidance of three functions) is equivalent
to finding a function with prescribed zeros, ones and poles
in the unit disc, counting multiplicity,
and prescribed jets at finitely many points. 

Interestingly, Nevanlinna \cite{Nevanlinna}
proposed a similar problem for meromorphic
functions in $\C$: to find necessary and sufficient conditions that
zeros, poles and $1$-points of a meromorphic functions must satisfy.
Some necessary conditions are known \cite{Nevanlinna,Winkler},
see also \cite{RY,O}. Most of these results are for meromorphic 
functions in the plane.

Consider the following examples.
\vspace{.1in}

1. (Blondel \cite{Blondel}) For which $\delta$ the following three transfer functions
in the unit disc are simultaneously stabilizable:
$$p_1=z^2/(z-\delta),\quad p_2(z)=z^2/(z+\delta),\quad p_3(z)=0\quad ?$$
The stabilizer $c$ has to be a rational function without poles
in the unit disc avoiding $1/p_i,\; i=1,2.$
This means $g(z)=z-c(z)z^2$ has to satisfy
$$g(0)=0,\quad g'(0)=1,\quad g(z)\neq \pm\delta,\quad |z|<1.$$
According to a result of Bermant \cite{Bermant} 
this is possible if and only if 
$$\delta\geq\delta_0:= 8\pi^2/\Gamma^4(1/4),$$
and the extremal function is not rational. 
This inequality gives a necessary condition
of simultaneous stabilizability.
Then an easy approximation argument shows
that $p_1,p_2$ and $p_3$ are simultaneously stabilizable if and only
if $\delta>8\pi^2/\Gamma^4(1/4).$
\vspace{.1in}

2. (Patel \cite{Patel}) For which $a>0$, the following three transfer functions
in the unit disc are simultaneously stabilizable:
$$p_1(z)=z,\quad p_2(z)=z^2/(z-a),\quad p_3(z)=0\quad ?$$
The stabilizer $c$ has to be a rational function without
poles in the unit disc satisfying
$$c(z)\neq 1/z,\quad c(z)\neq (z-a)/z^2,\quad |z|<1.$$
Introducing $g=(z-c(z)z^2)/a$ we rewrite this in the equivalent form:
\begin{equation}\label{hu}
g(z)\neq 1, \quad g(0)=0\leftrightarrow z=0, \quad g'(0)=1/a.
\end{equation}
The answer follows from a theorem of Caratheodory
\cite{Cara}, \cite{Nehari2}, \cite{Goluzin2}
\vspace{.1in}

{\em If a holomorphic function $g$ in the unit disc satisfies
(\ref{hu}) then 
$$|a|\geq 1/16.$$ 
There is a real holomorphic function for which
equality holds. This extremal function is not rational.}
\vspace{.1in}

A similar result, but with a smaller constant, was obtained
for the first time by Hurwitz \cite{Hurw}.

Now a simple approximation argument
shows that
the above three systems are simultaneously stabilizable
if and only if $a>1/16.$ 
This answers a question stated in \cite{Patel}.

Suppose that we wish to stabilize three transfer functions,
one of which avoids another. The problem is equivalent
to finding a rational function without zeros and poles in the unit disc,
which avoids one rational function $p$.
Such $c$ is
called a {\em bistable controller}. I am not aware of any practical
application of this ``bistability property''
by itself, but the desire to
control three systems with a single controller
is reasonable as explained
above.
The problem now is to find necessary and sufficient conditions
on a rational function $p$ for the existence of $c$ satisfying
(\ref{3}) and having no zeros and no poles in $D$.
Blondel \cite{Blondel} calls this 
``one of the major unsolved
problems of control theory''.

A special case of Theorem 1 above,
previously established by Blondel, says
that the problem in equivalent to: 
%

{\em Finding a rational function $w$
without poles in $D$, so that $1$-points of $w$ in $D$ 
and zeros of $w$ in $D$ are prescribed (with multiplicities).}
\vspace{.1in}

We refer to \cite{Blondel2,Blondel3} and the references in \cite{Blondel}
 for some necessary conditions
that zeros and $1$-points must satisfy. 

Only one universal restriction (independent
of degree) which zeros, poles and $1$-points of a rational function
must satisfy is known. It was found by Goldberg \cite{Gold}
and later independently by Blondel.

To state Goldberg's result, we introduce some notation.
Let $F_0$ be the class of all holomorphic functions $f$ in the rings
$$\rho(f)<|z|<1,$$
with the properties that $f(z)\not\in\{0,1,\infty\}$, and 
the indices (winding numbers) of the curve 
\begin{equation}\label{gamma}
\gamma(f)=\{ f(z):|z|=(1+\rho(f))/2\}
\end{equation}
about $0$ and $1$
are non-zero and distinct.
Let $F_4,F_3,F_2,F_1$ be the subsets of $F_0$
which consist of polynomials, rational,
holomorphic, and meromorphic functions in $D$, respectively, having
finite pairwise distinct numbers of zeros, poles and $1$-points.
We have $F_4\subset F_2\subset F_1\subset F_0$ and
$F_2\subset F_3\subset F_1 \subset F_0$. The constants $\rho(f)$ are
defined for $f\in F_j,\; 1\leq j\leq 4$ as
$$\rho(f)=\max\{|z|:f(z)\in\{0,1,\infty\}\}.$$
Now we put
$$A_j=\inf\{\rho(f):f\in F_j\},\quad 0\leq j\leq 4.$$
Evidently $A_0\leq A_1\leq A_3\leq A_4$ and $A_0\leq A_1\leq A_2\leq A_4$.
Goldberg's theorem says that
$$0<A_0=A_1=A_3<A_2=A_4.$$
Moreover, extremal functions exist for $A_0$ and $A_2$ but do not exist
for $A_1,A_3$ and $A_4$.

This result shows that if a holomorphic function in the unit disc
has finite, non-zero, distinct numbers of zeros and $1$-points,
then these zeros and one points cannot lie very close together.

So we have two absolute constants $0<A_0<A_2$ which are called
Goldberg's constants. The exact value of $A_0$ is known:
$$A_0=\exp(-\pi^2/(\log(3+2\sqrt{2}))\approx 0.003701599,$$
and for $A_2$ there are estimates
$$0.00587465<A_2\leq \mu\approx 0.0252896.$$
The constant $\mu$ and a function which corresponds to it are
conjectured to
be extremal for $A_2$; this function $h$ is described in detail
in \cite{BE}, and we
will give a short description below.

If the indices of the curve $\gamma$ about $0$ and $1$
are prescribed to be $N_0,N_1$,
we obtain constants $A_0(N_0,N_1)$.
One can obtain an exact value of $A_0(N_0,N_1)$, for any given $N_1>N_0>0$,
see \cite{BE}.

Being unable to prove that $A_2=\mu$, the authors of \cite{BE}
showed that $\mu$ is the solution of a restricted extremal problem:
\vspace{.1in}

\noindent
{\bf Theorem.} {\em A necessary and sufficient condition for the existence
of a holomorphic function $f$ in the unit disc, having no poles,
a single simple zero at $a$ and a single multiple $1$-point at $-a$
is that $|a|\geq\mu$. If $a=\mu$ this function is unique and transcendental.
If $|a|>\mu$ there exists a polynomial $f$ with the stated properties.}
\vspace{.1in}

Thus in the simplest case of one simple zero, one
multiple $1$-point and no poles, we have a necessary and sufficient condition
for the existence of a rational function with prescribed zeros, $1$-points
and poles in the unit disc.

This can be restated as a necessary and sufficient condition for
a stabilization problem as follows:
\vspace{.1in}

{\em The transfer function 
$$\frac{(z+a)^2}{z-a}$$
can be stabilized by a bistable controller if and only if $|a|>\mu$.}
\vspace{.1in}

In more complicated cases, there is no hope for such simple conditions. 
For example, Blondel \cite{Blondel} states the following 
problem\footnote{He even offered
a prize of $1$ kg of fine Belgian chocolate for this problem.
Nevertheless it is still wide open.}:
\vspace{.1in}

{\em For which $\delta>0$ there exists a rational
function which in the unit disc has no poles, a single simple zero at $0$,
and exactly two simple
$1$-points $\pm i\delta$?}
\vspace{.1in}

It is known that there exists $\delta_0>0$, 
with the property that such function
exists for $\delta\geq \delta_0$ and does not exist for $\delta<\delta_0$.

Evidently $\delta_0\geq A_2$ and it is not difficult to show that
this inequality is strict. 
The current world record \cite{Chang} for the estimate from above seems to be
$\delta_0<0.1148.$
The best known lower estimate is $0.01450779$. It can be obtained from
the estimate in \cite{Hempel} of the minimal length of a closed
hyperbolic geodesic in a twice punctured disc \cite{BE}.

In conclusion, we sketch the definition of the function
which is conjectured to be extremal for $A_2$.
The fundamental group $\Gamma$ of $\C\backslash\{0,1\}$
is a free group generated by
simple loops $A$ and $B$ around $0$ and $1$. Let $\Gamma^\prime$
be the subgroup generated by $A$ and $B^2$. It is also a free group
on two generators. Let $g:X\to\C\backslash\{0,1\}$ be the covering map
corresponding to this subgroup $\Gamma^\prime$, so that $\Gamma^\prime$
is the fundamental group of $X$. One can show that $X$ is a Riemann surface
which is conformally equivalent to the twice punctured disc, and we can identify
it with $D\backslash\{-\mu,\mu\}$ for some $\mu\in D$. Then $g$ becomes
a holomorphic function in $D$ which has one simple zero, say at $-\mu$
and one double $1$-point at $\mu$. We conjecture that $A_2=\mu.$
One can express our function $g$ in terms of solutions of a Lam\'e equation
and modular functions.

{\em Department of Mathematics, Purdue University

West Lafayette IN 47907

eremenko@math.purdue.edu}
\end{document}